\documentclass[10pt]{article}  
\usepackage[all]{xy}
\usepackage[mathscr]{euscript}
\usepackage{amsmath, amsfonts, amssymb, mathrsfs}
\usepackage{appendix}
\usepackage{amsthm}
\usepackage{comment}
\usepackage{textcomp}
\usepackage{setspace}
\usepackage{fix-cm}

 \usepackage{mathtools}
 \usepackage{nccmath}

\usepackage[alphabetic,lite]{amsrefs}

\usepackage{tikz-cd}

\usetikzlibrary{arrows,positioning,decorations.pathmorphing,
  decorations.markings
 }

\tikzset{inner sep=0pt, 
  root/.style={circle,draw,minimum size=7pt,thick}, 
  fatroot/.style={circle,draw,minimum size=10pt,thick}, 
  short root/.style={circle,fill,minimum size=7pt}, 
  doublearrow/.style={postaction={decorate}, 
  decoration={markings,mark=at position .7
  with {\arrow{angle 60}}},double distance=3pt,thick}
}

\newcommand{\bC}{\mathbb C}

\newcommand{\bZ}{\mathbb Z}

\newcommand{\bQ}{\mathbb Q}
\newcommand{\bR}{\mathbb R}

\newcommand{\cB}{\mathcal B}
\newcommand{\cO}{\mathcal O}

\newcommand{\frakg}{\mathfrak g}

\newcommand{\la}{\langle}
\newcommand{\ra}{\rangle}

\newcommand{\al}{\alpha}
\newcommand{\be}{\beta}

\newcommand{\ep}{\epsilon}
\newcommand{\lam}{\lambda}
\newcommand{\varep}{\varepsilon}

\DeclareMathOperator{\Ad}{Ad}

\DeclareMathOperator{\Hom}{Hom}
\DeclareMathOperator{\Sym}{Sym}

\DeclareMathOperator{\GL}{GL}

\DeclareMathOperator{\SL}{SL}

\DeclareMathOperator{\SO}{SO}

\DeclareMathOperator{\Irr}{Irr}

\DeclareMathOperator{\ind}{ind}

\DeclareMathOperator{\Stab}{Stab}

\newcommand{\cA}{\mathcal A}

\newcommand{\frakf}{\mathfrak f}

\newcommand{\Gk}{{G}} 
\newcommand{\Gx}{{\mathsf{G}_x}} 
\newcommand{\Vx}{\mathsf{V}_x}
\newcommand{\checkVx}{\check{\mathsf{V}}_x}

\newcommand{\resk}{\mathfrak{f}}
\newcommand{\resK}{\mathfrak{F}}

\DeclareMathOperator{\Frob}{Fr} 
\newcommand{\Weil}{\mathcal{W}} 
\newcommand{\Inert}{\mathcal{I}} 
 
\DeclareMathOperator{\Frobk}{Fr_0} 
\DeclareMathOperator{\LL}{LLC_m} 
\DeclareMathOperator{\LLk}{LLC_1} 
\DeclareMathOperator{\LLG}{LLC_G}

\newcommand{\GZ}{\mathscr{G}}
\newcommand{\TZ}{\mathscr{T}}
\newcommand{\VZ}{\mathscr{V}}
\newcommand{\MZ}{\mathscr{M}}

\DeclareMathOperator{\im}{im}
\DeclareMathOperator{\Lie}{Lie}
\DeclareMathOperator{\spn}{span}
\DeclareMathOperator{\Ht}{ht}
\DeclareMathOperator{\Deg}{Deg} 
\DeclareMathOperator{\Rep}{Rep} 

\newcommand{\princ}{\varphi_0} 

\theoremstyle{plain}
\input xy
\xyoption{all}
\newtheorem{Thm}{Theorem}[section]

\newtheorem{Lem}[Thm]{Lemma}
\newtheorem{Prop}[Thm]{Proposition}
\newtheorem{Cor}[Thm]{Corollary}

\theoremstyle{remark}
\newtheorem{remark}[Thm]{Remark}
\numberwithin{equation}{section}
\numberwithin{paragraph}{section}

\AtEndDocument{\bigskip{\footnotesize 
\textsc{Beth Romano}~~ \texttt{beth.romano@kcl.ac.uk}\\\textsc{Department of Mathematics, King's College London,
WC2R 2LS, UK}}}

\begin{document}

\author{Beth Romano}
\title{On input and Langlands parameters for epipelagic representations}
\date{}
\maketitle

\thispagestyle{empty}

\begin{abstract}
A paper of Reeder--Yu gives a construction of epipelagic supercuspidal representations of $p$-adic groups.
The input for this construction is a pair $(\lam, \chi)$ where $\lam$ is a stable vector in a certain representation coming from a Moy--Prasad filtration, and $\chi$ is a character of the additive group of the residue field. We say two such pairs are equivalent if the resulting supercuspidal representations are isomorphic. In this paper we describe the equivalence classes of such pairs.
As an application, we give a classification of the simple supercuspidal representations for split adjoint groups. Finally, under an assumption about unramified base change, we describe properties of the Langlands parameters associated to these simple supercuspidals, showing that they have trivial L-functions and minimal Swan conductors, and showing that each of these simple supercuspidals lies in a singleton L-packet.
\end{abstract}

\tableofcontents

\section{Introduction}

In \cite{ReederYu}, Reeder and Yu 
defined the notion of \textit{epipelagic} supercuspidal representations of $p$-adic groups and gave a construction of a subset of these representations. 
Epipelagic representations have small positive depth, and so they give a useful test case for questions about the positive-depth part of the local Langlands correspondence (for example questions about preservation of depth, cf. \cite{ABPSDepth}, \cite{AMPS}). 
The construction of \cite{ReederYu} is particularly useful because it works in a uniform way for all residue characteristic and all semisimple groups, including exceptional groups. 
It also gives a way to build supercuspidal representations that cannot be constructed via the methods of \cite{Yu} or \cite{Stevens}. Since the construction of \cite{Yu} is exhaustive for residue characteristic $p$ not dividing the order of the Weyl group (\cite{Fintzen}), and the construction of \cite{Stevens} is exhaustive for classical groups when $p \neq 2$, the case of small residue characteristic for exceptional groups can be seen as the final frontier in the classification of supercuspidals. 

Yet \cite{ReederYu} doesn't give information about when two distinct inputs for their construction yield the same representations. This information is needed in order to parametrize the resulting representions, which is an important step toward building an explicit local Langlands correspondence in cases where it is not already known. The main result of this paper fills this gap, giving a criterion for when distinct inputs for the construction of \cite{ReederYu} yield isomorphic supercuspidal representations (Theorem \ref{thm-input}). 

As an application of this result, I give an explicit construction of all \textit{simple supercuspidal representations} for split adjoint groups (Proposition \ref{prop-decomp-pix} and Corollary \ref{cor-formal-deg}). 
Simple supercuspidals 
form a special case of epipelagic representations: they 
were first introduced in \cite{GrossReeder} for simply connected groups, and the construction of \cite{ReederYu} can be see as a generalization of this earlier construction. 
Simple supercuspials have inspired a wealth of research, for example \cite{KalethaSimple}, \cite{Adrian}, \cite{AdrianKaplan1}, \cite{Oi},\cite{HenniartOi}, but they have not been fully explored in the case of adjoint groups, and most work on them has been restricted to classical groups and $\GL_n$.

To build a local Langlands correspondence for a reductive group $\Gk$ defined over a nonarchimedean local field $k$, one must associate a Langlands parameter to each smooth irreducible representation of $\Gk(k)$. This has been acheived in many cases for classical groups, notably for $\GL_n$ (\cite{HarrisTaylor}, \cite{Henniart}), for inner forms of $\SL_n$ (\cite{HiragaSaito}, \cite{ABPS}) and for symplectic and special orthogonal groups (\cite{Arthur}). Recently 
the local Langlands correspondence has been constructed for the group $G_2$ (\cite{AubertXu} for residue characteristic $p \neq 2, 3$; \cite{GanSavin} for all residue characteristic), but it is not fully understood for any other exceptional group. Again, the construction of \cite{ReederYu} gives a way to begin to fill the gaps in our understanding, particularly for exceptional groups over fields with small residue characteristic.

In the case when $\Gk$ is simply connected, \cite[Section 7]{ReederYu} describes properties of Langlands parameters that should correspond to the representations of \cite{ReederYu}.
Yet Langlands parameters for the representations of \cite{ReederYu} have only been explored beyond the simply connected case
in a few examples (cf. \cite{Adrian}, \cite{AdrianKaplan2}, \cite{Oi19}, \cite{Oi}); these examples are 
restricted to classical groups and make assumptions on residue characteristic.
The second goal of this paper is to describe properties of the Langlands parameters associated to simple supercuspidal representations for split adjoint $\Gk$. In particular, we show that these Langlands parameters have trivial adjoint L-function, trivial image of $\SL_2(\bC)$, and minimal Swan conductor (Theorem \ref{thm-LP} and Corollary  \ref{cor-swan}). We also prove that in contrast to the simply connected case, the L-packets are singletons for simple supercuspidal representations of split adjoint groups (Remark \ref{rem-packet}).
Some of these results extend those of \cite{Adrian} and \cite{Oi} for $\SO_{2n + 1}$. We emphasize that our proofs make no assumptions about $k$ and are uniform for all split adjoint groups, but to use the results of \cite{GrossReeder}, we make an assumption about unramified base change (see Section \ref{section-LLC}).

We now describe the results of this paper in more detail. Let $k$ be a nonarchimedean local field with residue field $\resk$, and let $\Gk$ be a connected semisimple algebraic group over $k$ that splits over a tamely ramified extension of $k$. Given a rational point $x$ in the apartment associated to a maximal $k$-split torus in $\Gk$, we may consider the corresponding Moy--Prasad filtration
\begin{equation*}
\Gk(k)_{x, 0} \gneq \Gk(k)_{x, r_1} \gneq \Gk(k)_{x, r_1} \gneq \dots
\end{equation*}
where $(r_1, r_2, \dots)$ is an increasing sequence of positive real numbers. The first quotient
$\Gx(\resk) = \Gk(k)_{x, 0}/\Gk(k)_{x, r_1}$  forms the $\resk$-points of a reductive group over $\resk$, and the quotient $\Vx(\resk) = \Gk(k)_{x, r_1}/\Gk(k)_{x, r_2}$ is an $\resk$-vector space with a natural action of $\Gx(\resk)$ induced by conjugation in $\Gk$. The input for the construction of \cite{ReederYu} is a pair $(\lam, \chi)$, where $\lam$ is a stable vector (in the sense of geometric invariant theory, see Section \ref{section-review-epi} for details) in the dual representation $\checkVx(\resk)$ and $\chi: \resk^+ \to \bC^\times$ is a nontrivial character of the additive group $\resk^+$. 

The compactly induced representation
\begin{equation*}
	\pi_x(\lam, \chi) := \ind_{\Gk(k)_{x, r_1}}^{G(k)} (\chi \circ \lam) 
	\end{equation*}
decomposes into finitely many irreducible supercuspidal representations (\cite[Proposition 2.4]{ReederYu}), which are the epipelagic representations of interest in this paper.
Write $I(x)$ for the set of input pairs $(\lam, \chi)$ as described above and $\Irr(\lam, \chi)$ for the set of irreducible subrepresentations of $\pi_x(\lam, \chi)$, considered up to isomorphism. Let $H_x$ be the stabilizer of $x$ in $\Gk(k)$. 

\begin{Thm}[Theorem \ref{thm-input}]\label{thm-input-intro}
Suppose $(\lam, \chi), (\lam', \chi') \in I(x)$. Choose $c \in \resk^\times$ such that $\chi' = c\chi$.
\begin{enumerate}
\item[1.] The sets $\Irr(\lam, \chi)$ and $\Irr(\lam', \chi')$ are either disjoint or equal. If they are equal, then $\pi_x(\lam, \chi) \simeq \pi_x(\lam', \chi')$.
\item[2.] 
We have $\pi_x(\lam, \chi) \simeq \pi_x(\lam', \chi')$ if and only if there exists $g \in H_x$ such that $g \cdot \lam = c\lam'$. In particular, $\pi_x(\lam, \chi) \simeq \pi_x(\lam', \chi)$ if and only if $\lam$ and $\lam'$ are in the same $H_x$-orbit.
\end{enumerate}
\end{Thm}

Now assume $\Gk$ is split, absolutely simple, and adjoint, and assume that $x$ is the barycenter of an alcove. In this case we define an invariant polynomial $\Delta$ on $\checkVx(\resk)$, and we show that two vectors $\lam, \lam' \in \checkVx(\resk)$ are in the same $H_x$-orbit if and only if $\Delta(\lam) = \Delta(\lam')$ (Proposition \ref{prop-orbit-reps}). Using this, we prove the following result. Here $\Omega$ is the stabilizer of an alcove in the extend affine Weyl group of $\Gk$, $q$ is the cardinality of $\resk$, and $J_x = \Gk(k)_{x, r_1}$ is the pro-unipotent radical of the Iwahori subgroup $\Gk(k)_{x, 0}$.

\begin{Cor}[Corollary \ref{cor-formal-deg}]\label{cor-ssc-intro}
Suppose $\Gk$ is split, absolutely simple, and adjoint, and that $x$ is the barycenter of an alcove. 
Then the construction of \cite{ReederYu} gives a total of $(q - 1)\lvert \Omega \rvert$ supercuspidal representations. With respect to a Haar measure $\mu$ on $\Gk(k)$, each of these representations has 
 formal degree $\frac{1}{\lvert \Omega \rvert \mu(J_x)}$. 
 \end{Cor}
 
 The representations are described in more detail in Proposition \ref{prop-decomp-pix}.
 
 Now we turn to Langlands parameters. With $\Gk$ split and adjoint as above, let $\ell$ be the rank of $\Gk$. Let $G^\vee$ be the ($\bC$-points of the) dual group of $\Gk$, let $Z^\vee$ be the center of $G^\vee$, and let $\frakg^\vee = \Lie(G^\vee)$. Let $\pi$ be one of the supercuspidal representations described in Corollary \ref{cor-ssc-intro}. 
 Here $\al(\varphi)$ is the Weil--Deligne version of the Artin conductor (see Section \ref{section-LP-defs}).

\begin{Thm}[Theorem \ref{thm-LP}]\label{thm-LP-intro}
Suppose $(\varphi, \rho)$ is the enhanced Langlands parameter corresponding to $\pi$ under the local Langlands correspondence. Then
 \begin{enumerate}
\item[1.] 
$\al(\varphi) = \dim \frakg^\vee + \ell$. 
\item[2.]  $\varphi(\SL_2(\bC)) = 1$.
\item[3.] The adjoint L-function $L(\varphi, s)$ is trivial. 
\item[4.] $Z_{G^\vee}(\im \varphi) = Z^\vee$, and $\rho$ is the trivial representation.
\end{enumerate}
\end{Thm}

As a corollary we see that $\pi$ lies in a singleton L-packet (Remark \ref{rem-packet}) and that $\varphi$ has minimal Swan conductor (Corollary \ref{cor-swan}). 

Note that to prove Theorem \ref{thm-LP-intro} we assume the local Langlands correspondences satisfies two properties. The first is the formal degree conjecture of \cite{HiragaIchinoIkeda}. (A proof of this conjecture has been announced for classical groups by Beuzart-Plessis.)
To explain the second property, let $k'/k$ be a finite unramified extension. 
We show in Section \ref{section-unram} that there is a way to relate to every simple supercuspidal representation $\pi$ of $\Gk(k)$ a corresponding simple supercuspidal $\pi'$ of $\Gk(k')$. 
For certain extensions $k'/k$, we assume we can obtain the Langlands parameter for $\pi'$ by restriction from the Langlands parameter of $\pi$ (see Section \ref{section-LLC}). This assumption is analogous to one made in \cite{GrossReeder} when considering simple supercuspidals for simply connected groups, and while it is not explicitly stated, a similar assumption seems necessary for the proofs of \cite[Section 7]{ReederYu}.

\begin{remark}
Not long after this paper was posted online, the paper \cite{AHKO}, which focuses on simple supercuspidal representations for classical groups, was posted. There is some overlap in the results of these two papers. In particular, \cite{AHKO} gives alternate proofs for some of the results of sections \ref{section-ssc} and \ref{section-LPs} in the special case when $\Gk = \SO_{2n + 1}$.
\end{remark}

\subsection{Structure of the paper}

In Section \ref{section-review-epi} we review the construction of \cite{ReederYu}, and in Section \ref{section-input-result} we give a criterion for two inputs for this construction to produce isomorphic supercuspidal representations. 

In Section \ref{section-ssc}, we consider the case when $\Gk$ is a simple, split, adjoint group, and $x$ is the barycenter of an alcove in the apartment associated to a maximal split torus of $\Gk$. In order to use the results of Section \ref{section-input}, in Section \ref{section-Hx} we completely describe the orbits for the action of $H_x$ on $\checkVx(\resk)$ in terms of an explicit invariant polynomial. 
We consider $(\chi, \lam), (\chi', \lam') \in I(x)$ to be equivalent if $\pi_x(\lam, \chi) \simeq \pi_x(\lam', \chi')$. In Propositions \ref{prop-orbit-reps} and \ref{prop-ssc-classification} we completely describe the equivalence classes for $I(x)$ and give representatives for the equivalence classes. In Section \ref{section-decomp-ssc} we decompose the representations $\pi_x(\lam, \chi)$. Corollary \ref{cor-formal-deg} gives the number of representations resulting from the construction of \cite{ReederYu} in this case. 

In Section \ref{section-unram} we consider finite unramified extensions $k'/k$. We show that there is a natural way to relate simple supercuspidal representations of $\Gk(k)$ to certain simple supercuspidal representations of $\Gk(k')$. This relation will allow us to relate Langlands parameters for the corresponding representations in an analogous way to \cite{GrossReeder}.

In Section \ref{section-LPs} we move on to Langlands parameters. In Section \ref{section-LP-defs} we review some background about Langlands parameters and L-functions. In Section \ref{section-LLC}, we show that if $\pi$ is one of the irreducible supercuspidal representations constructed in Section \ref{section-ssc}, then the corresponding Langlands parameter has trivial adjoint L-function, trivial image of $\SL_2(\bC)$, and minimal Swan conductor.

\subsection{Notation}

Suppose $G$ is a group and $H$ is a subgroup of $G$. Given a function $\chi$ defined on $H$ and an element $g \in G$, we write $\chi^g$ for the function on $H^g := g^{-1}Hg$ defined by $\chi^g(x) = \chi(gxg^{-1})$ for all $x \in H^g$. In other words, $\chi^g(g^{-1}hg) = \chi(h)$ for all $h \in H$. Similarly, if $(\rho, V)$ is a representation of $H$, then $\rho^g$ defines a representation of $H^g$ on the same vector space $V$. We sometimes write $V^g$ for this representation. 

If $G$ is a locally profinite group, $H$ is a closed subgroup of $G$, and $\tau$ is a smooth representation of $H$, we write $\ind_H^G\tau$ for the compact induction of $\tau$.

Throughout the paper, we fix a nonarchimedean local field $k$ with finite residue field $\resk$ of cardinality $q$ and characteristic $p$. We write $K$ for the maximal unramified extension of $k$ in a fixed separable closure. Let $\resK$ be the residue field of $K$, so $\resK$ is an algebraic closure of $\resk$.
Let $\Gk$ be a connected semisimple algebraic group over $k$. We assume that $\Gk$ splits over a tamely ramified extension of $k$. Let $\Frobk$ be the Frobenius automorphism on $G(K)$ arising from the given $k$-structure on $\Gk$: we have $G(K)^{\Frobk} = G(k)$.

Choose a maximal $k$-split torus $S$ in $G$, and let $T$ be a maximal $K$-split torus of $G$ containing $S$. We choose $T$ to be defined over $k$. Let $\cB(G, K)$ be the Bruhat--Tits building of $G(K)$, and let $\cA(T, K) \subset \cB(G, K)$ be the apartment corresponding to $T$. 
Let $\Psi$ denote the affine $K$-roots of $T$ in $\Gk$.
We may identify the apartment $\mathcal{A}(S, k)$ corresponding to $S$ in the Bruhat--Tits building of $\Gk(k)$ with $\mathcal{A}(T, K)^{\Frobk}$. A point $x \in \cA(S, k)$ is called \textit{rational} if $\psi(x) \in \bQ$ for every $\psi \in \Psi$.

\subsection*{Acknowledgements}

I'd like to thank Guy Henniart, Mark Reeder, and Maarten Solleveld for helpful comments on a draft of this paper, as well as the anonymous referee for their thoughtful feedback.

\section{Input in the Reeder--Yu construction}\label{section-input}

\subsection{Review of the construction}\label{section-review-epi}

We now briefly review the construction of supercuspidal representations in \cite[Section 2]{ReederYu}. 
Fix a rational point $x \in \mathcal{A}(S, k)$.
For $r \geq 0$ let $\Gk(K)_{x, r}$ be the Moy--Prasad filtration subgroup as defined in \cite{MP1, MP2}.
The groups $\Gk(K)_{x, r}$ are normal in $\Gk(K)_{x, 0}$ for all $r \geq 0$. 
 Set 
\begin{equation*}
G(K)_{x, r+} = \underset{s > r}\bigcup G(K)_{x, s}.
\end{equation*} 
The groups $G(K)_{x, r}$ are stable under the action of $\Frobk$,
and we define subgroups $G(k)_{x, r}$ of $G(k)$ by $G(k)_{x, r} = (G(K)_{x, r})^{\Frobk}$ (and similarly for $G(k)_{x, r+}$). 

Let $r(x)$ be the smallest positive value in the set $\{ \psi(x) \mid \psi \in \Psi \}$. Let
\begin{align*}
\Gx(\resK) &= \Gk(K)_{x, 0}/G(K)_{x, r(x)},\\
\Vx(\resK) &= \Gk(K)_{x, r(x)}/G(K)_{x, r(x)+}.
\end{align*}
Then $\Gx(\resK)$ forms the $\resK$-points of a reductive group $\Gx$ over $\resK$, and $\Vx(\resK)$ is a finite-dimensional vector space over $\resK$.
The action of $\Gk(K)_{x, 0}$ by conjugation induces an action of
$\Gx(\resK)$ on $\Vx(\resK)$. Let $\checkVx(\resK)$ be the dual representation. A vector $\lam \in \checkVx(\resK)$ is called \textit{stable} (in the sense of geometric invariant theory, cf. \cite{Mumford}) if its orbit $\Gx(\resK)\cdot \lam$ is closed and its stabilizer in $\Gx(\resK)$ is finite.

Let 
\begin{align*}
J_x &= G(k)_{x, r(x)},\\
J_x^+ &= G(k)_{x, r(x)+}.
\end{align*} 
The action of $\Frobk$ on $\Gk(K)$ induces actions on the quotients $\Gx(\resK)$ and $\Vx(\resK)$. 
Thus we may think of $\Gx$ as a group defined over $\resk$, and of $\Vx$ as an algebraic representation of $\Gx$ over $\resk$.
We say that a vector $\lam \in \checkVx(\resk)$ is stable if it is stable as an element of $\checkVx(\resK) = \checkVx(\resk) \otimes_{\resk} \resK$  in the sense defined above.

Given a stable vector  $\lam \in \checkVx(\resk)$ and a nontrivial character $\chi: \resk^+ \to \bC^\times$, we consider the composition 
\begin{equation*}
\chi_\lam := \chi \circ \lam: \Vx \to \bC^\times
\end{equation*}
 as a character of $J_x$ that is trivial on $J_x^+$. Then the compactly induced representation
	\begin{equation*}
	\pi_x(\lam, \chi) := \ind_{J_x}^{G(k)} (\chi_\lam) 
	\end{equation*}
	is a direct sum of finitely many irreducible supercuspidal representations of $\Gk(k)$ of depth $r(x)$ (\cite[Propostion 2.4]{ReederYu}). These are the epipelagic supercuspidal representations of interest in this paper.

Proposition 2.4 in \cite{ReederYu} also describes the irreducible summands of $\pi_x(\lam, \chi)$ more concretely. 
Consider the stablizer $\Stab_{\Gk(K)}(x)$ for the action of $\Gk(K)$ on $\cB(G, K)$. This stabilizer contains the parahoric subgroup $\Gk(K)_{x, 0}$ with finite index.
The actions of $\Gk(K)_{x, 0}$ on $\Vx(\resK)$ and $\checkVx(\resK)$ extend to actions of $\Stab_{\Gk(K)}(x)$.
 Since $\Stab_{\Gk(K)}(x)/\Gk(K)_{x, 0}$ is finite, $\Stab_{\Gk(K)}(x)$ preserves the set of stable vectors in $\checkVx(\resK)$.

Let $H_x = (\Stab_{\Gk(K)}(x))^{\Frobk}$. 
Given a stable vector $\lam \in \checkVx(\resk)$, let $H_\lam = \Stab_{H_x}(\lam)$.
Then the irreducible subrepresentations of $\pi_x(\lam, \chi)$ are the representations of the form $\ind_{H_\lam}^{\Gk(k)}(\rho)$, where $\rho$ is an irreducible subrepresentation of $\ind_{J_x}^{H_\lam}\chi_\lam$.

\subsection{On choice of input}\label{section-input-result}

Given a rational point $x \in \mathcal{A}(S, k)$ as above, let $I(x)$ be the set of pairs $(\lam, \chi)$ such that $\lam \in \checkVx(\frakf)$ is stable and $\chi: \mathfrak{f}^+ \to \bC^\times$ is a nontrivial character. Thus $I(x)$ is the set of possible inputs for the construction of \cite{ReederYu}. Given $(\lam, \chi) \in I(x)$, let $\Irr(\lam, \chi)$ be the set of irreducible supercuspidal representations (considered up to isomorphism) of $\Gk(k)$ appearing as summands of $\pi_x(\lam, \chi)$ as defined above. 
We note for future use that if $g \in H_x$ and $(\lam, \chi) \in I(x)$, then $\chi_{g\cdot \lam} = (\chi_\lam)^{g^{-1}}$.

\begin{Lem}\label{lem-twist}
Let $g \in \Gk(k)$, and suppose $(\lam, \chi), (\lam', \chi') \in I(x)$. If $\chi_\lam = (\chi'_{\lam'})^g$ on $J_x \cap J_x^g$, then $g \in H_x$. 
\end{Lem}

\textit{Proof}. We follow a similar logic to the proof of \cite[Proposition 2.4]{ReederYu}. 
First suppose $g \in N_{\Gk(k)}(S)$, which implies that $g$ preserves $\mathcal{A}(S, k)$ under the action of $\Gk(k)$ on the Bruhat--Tits building of $\Gk$. 
We have that $(J_x^+)^g = J^+_{g^{-1}\cdot x}$.
 Since $\chi_{\lam'}'$ is trivial on $J_x^+$, we have that $\chi_\lam$ is trivial on $(J_x^+)^g \cap J_x$. 
Let $\mathsf{V}_{x, g^{-1}\cdot x}(\resk)$ be the image of $J_{g^{-1}\cdot x}^+ \cap J_x$ in $\Vx(\resk)$, and let $\mathsf{V}_{x, g^{-1}\cdot x}(\resK)$ be the image of 
\begin{equation*}
\Gk(K)_{g^{-1}\cdot x, r(x)+} \cap \Gk(K)_{x, r(x)}
\end{equation*}
 in $\Vx(\resK)$. 
Since $\mathsf{V}_{x, g^{-1}\cdot x}(\resk)$ contains a basis for $\mathsf{V}_{x, g^{-1}\cdot x}(\resK)$, the map
 $\chi \circ \lam$ is trivial on $\mathsf{V}_{x, g^{-1}\cdot x}(\resK)$. Since $\mathsf{V}_{x, g^{-1}\cdot x}(\resK)$ is a subspace of $\Vx(\resK)$,  this implies that $\lam$ vanishes on $\mathsf{V}_{x, g^{-1}\cdot x}(\resK)$. 
By \cite[Lemma 2.3]{ReederYu}, we have $g^{-1}\cdot x = x$, i.e. $g \in H_x$.

Now we consider an arbitrary element $g \in \Gk(k)$. Using the affine Bruhat decomposition \cite[3.3.1]{Tits} and the fact that $H_x$ contains an Iwahori subgroup, 
we may write $g$ as $anb$ where $a, b \in H_x$ and $n \in N_{\Gk(k)}(S)$. 
Suppose $y \in J_x \cap J_x^n$. Since $H_x$ normalizes $J_x$,
 we have that $b^{-1}yb \in J_x$. We claim that $b^{-1}yb \in J_x^g$. Indeed, we have
 \begin{equation*}
 b^{-1}yb = g^{-1}(anyn^{-1}a)g.
 \end{equation*} Since $y \in J_x^n$, we have $nyn^{-1} \in J_x$, and $a$ normalizes $J_x$, which implies the claim.
 
Using this, we have that
\begin{align*}
(\chi^{b^{-1}})(y) &= \chi(b^{-1}yb)\\
&= \chi'(anyn^{-1}a^{-1})\\
&= ((\chi')^a)^n(y)
\end{align*}
i.e. $\chi_\lam^{b^{-1}} = (\chi_{\lam'}^a)^n$ on $J_x \cap J_x^n$. Now $b\cdot \lam$ and $a^{-1}\cdot \lam' \in \checkVx(\resk)$ are stable vectors with $\chi_\lam^{b^{-1}} = \chi_{b\cdot \lam}$ and $\chi_{\lam'}^a = \chi_{a^{-1}\cdot \lam'}$. 
Since 
\begin{equation*}
\chi_{b\cdot \lam} = (\chi_{a^{-1}\cdot \lam'})^n
\end{equation*}
on $J_x \cap (J_x)^n$,
we have $n \in H_x$ by the first case proven above. This implies $g \in H_x$. \qed

To state the theorem, recall that $\resk^\times$ acts transitively on the set of nontrivial characters of $\resk^+$. Explicitly, fix a nontrivial character $\chi: \resk^+ \to \bC^\times$. For $c \in \resk^\times$, let $c\chi: \resk^+ \to \bC^\times$ be defined by $c\chi(y) = \chi(cy)$ for all $y \in \resk^+$. Then the set of nontrivial characters of $\resk^+$ is given by $\{c\chi \mid c \in \resk^\times\}$.

\begin{Thm}\label{thm-input}
Suppose $(\lam, \chi), (\lam', \chi') \in I(x)$. Choose $c \in \resk^\times$ such that $\chi' = c\chi$.
\begin{enumerate}
\item[1.] The sets $\Irr(\lam, \chi)$ and $\Irr(\lam', \chi')$ are either disjoint or equal. If they are equal, then $\pi_x(\lam, \chi) \simeq \pi_x(\lam', \chi')$.
\item[2.] 
We have $\pi_x(\lam, \chi) \simeq \pi_x(\lam', \chi')$ if and only if there exists $g \in H_x$ such that $g \cdot \lam = c\lam'$. In particular, $\pi_x(\lam, \chi) \simeq \pi_x(\lam', \chi)$ if and only if $\lam$ and $\lam'$ are in the same $H_x$-orbit.
\end{enumerate}
\end{Thm}

\textit{Proof}. 
First suppose there exists $g \in H_x$ such that $g \cdot \lam = c\lam'$. We have 
\begin{equation*}
\pi_x(\lam', c\chi) = \pi_x(c\lam', \chi) = \pi_x(g\cdot \lam, \chi) = \ind_{J_x}^{\Gk(k)}\chi_{g\lam} = \ind_{J_x}^{\Gk(k)}(\chi_\lam)^{g^{-1}}.
\end{equation*}
 By transitivity of induction, 
we have
\begin{equation*}
\ind_{J_x}^{\Gk(k)}(\chi_\lam)^{g^{-1}} = \ind_{H_x}^{\Gk(k)}\ind_{J_x}^{H_x}(\chi_\lam)^{g^{-1}}.
\end{equation*} 
Now $\ind_{J_x}^{H_x}(\chi_\lam)^{g^{-1}}$ is isomorphic to $(\ind_{J_x}^{H_x} \chi_\lam)^{g^{-1}}$ as an $H_x$-module (it's not hard to check that $f \mapsto f^{g^{-1}}$ defines an $H_x$-module isomorphism).
Since $g \in H_x$, we have $(\ind_{J_x}^{H_x}\chi_\lam)^{g^{-1}} \simeq \ind_{J_x}^{H_x}\chi_\lam$, and thus $\pi_x(g\cdot\lam, \chi) \simeq \pi_x(\lam, \chi)$.

Now suppose $\Irr(\lam, \chi) \cap \Irr(\lam', \chi')$ is nonempty. We must show that there exists $g \in H_x$ such that $g\cdot \lam = c\lam'$. 
Since $\Irr(\lam, \chi) \cap \Irr(\lam', \chi')$ is nonempty, there exists an irreducible $H_{\lam}$-subrepresentation $\rho$ of $\ind_{J_x}^{H_{\lam}}(\chi_{\lam})$ and an irreducible $H_{\lam'}$-subrepresentation $\rho'$ of $\ind_{J_x}^{H_{\lam'}}(\chi'_{\lam'})$ such that
\begin{equation*}
\Hom_{\Gk(k)}(\ind_{H_\lam}^{\Gk(k)} \rho, \ind_{H_{\lam'}}^{\Gk(k)} \rho') \neq 0.
\end{equation*} 
By \cite{KutzkoMackey}, this means there exists $g \in \Gk(k)$ such that 
\begin{equation*}
\Hom_{H_\lam \cap (H_{\lam'})^g}(\rho, (\rho')^g) \neq 0.
\end{equation*}
 As explained in \cite[Section 2.1]{ReederYu}, $\rho$ is $\chi_\lam$-isotypic when restricted to $J_x$, and $\rho'$ is $\chi'_{\lam'}$-isotypic. Thus $(\rho')^g$ is $(\chi'_{\lam'})^g$-isotypic when restricted to $(J_x)^g$. We conclude that $\chi_\lam = (\chi'_{\lam'})^g$ on $J_x \cap (J_x)^g$. By Lemma \ref{lem-twist}, we must have $g \in H_x$. 
Now 
\begin{equation*}
(\chi'_{\lam'})^g = \chi'_{g^{-1} \cdot \lam'} = (\chi_c)_{g^{-1} \cdot \lam'} = \chi_{g^{-1} \cdot (c\lam')},
\end{equation*} 
so $\chi_\lam = \chi_{g^{-1}\cdot (c\lam)}$. This tells us that $(\lam - (g^{-1}\cdot c\lam'))(y) \in \ker \chi$ for all $y \in J_x$. Since $\chi$ is nontrivial and any nontrivial element of $\checkVx(\resk)$ is surjective, we must have $\lam = g^{-1}\cdot (c\lam')$, i.e. $g\lam = c\lam'$ as claimed. \qed

The following corollary of the theorem is useful for applications, particularly for counting the number of supercuspidal representations resulting from the construction of \cite{ReederYu}. Given $(\lam, \chi), (\lam', \chi') \in I(x)$, we write $(\lam, \chi) \sim (\lam', \chi')$ if $\Irr(\lam, \chi) \cap \Irr(\lam', \chi')$ is nonempty, which by the theorem is equivalent to $\pi_x(\lam, \chi) \simeq \pi_x(\lam', \chi')$. It's clear that $\sim$ forms an equivalence relation on $I(x)$. 

\begin{Cor}\label{cor-equiv-classes}
\begin{enumerate}
\item[1.] Suppose  $(\lam, \chi), (\lam', \chi') \in I(x)$, and choose $c \in \resk^\times$ such that $\chi' = c\chi$. We have $(\lam, \chi) \sim (\lam', \chi')$ if and only if there exists $g \in H_x$ such that $g\cdot \lam = c\lam'$.
\item[2.] Fix a nontrivial character $\chi_0$ of $\resk^+$. Each $\sim$-equivalence class in $I(x)$ contains a representative of the form $(\lam, \chi_0)$.
\item[3.] The number of equivalence classes in $I(x)$ is equal to the number of $H_x$-orbits on the stable vectors in $\checkVx(\resk)$.
\end{enumerate}
\end{Cor}

\textit{Proof}. 
The first statement is a reformulation of Theorem \ref{thm-input}. For the second statement, suppose $(\lam, \chi) \in I(x)$, and choose $c \in \resk^\times$ such that $\chi = c\chi_0$. Then $\chi_0 = c^{-1}\chi$. We have $\lam =  c^{-1}(c\lam)$, so by 1, we have $(\lam, \chi) \sim (c\lam, \chi_0)$. The third statement follows easily from the first two. Indeed, fix a nontrivial character $\chi_0$ of $\resk^+$, and let $S$ be a set of stable-orbit representatives for the action of $H_x$ on $\checkVx(\resk)$. If $\lam \in S$ then clearly $(\lam, \chi_0) \in I(x)$. And if $(\lam, \chi_0) \sim (\lam', \chi_0)$, then by the first part of the lemma, $\lam$ and $\lam'$ are in the same $H_x$-orbit. This shows that $\{(\lam, \chi_0) \mid \lam \in S\}$ is a complete set of $\sim$-equivalence class representatives. \qed

\begin{remark}
Note that if $\Gk$ is simply connected, then $H_x = \Gk_{x, 0}(k)$, and so the number of equivalence classes in $I(x)$ is equal to the number of stable orbits for $\Gx(\resk)$ acting on $\checkVx(\resk)$. But in general the orbits for $H_x$ in $\checkVx(\resk)$ are not necessarily the same as the orbits for $\Gx(\resk)$.
\end{remark}

\section{Simple supercuspidal representations for adjoint groups}\label{section-ssc}

In this section, we consider the supercuspidal representations coming from the construction of \cite{ReederYu} in the case when $G$ is split and $x$ is the barycenter of an alcove in $\cA(S, k)$. When $\Gk$ is simply connected, these representations were constructed in \cite{GrossReeder}, where they were named \textit{simple supercuspidal representations}. This name now more commonly refers to any representations coming from the \cite{ReederYu} construction for this choice of $x$. 
 Here we restrict ourselves to the case when $\Gk$ is adjoint. The case when $\Gk$ is simply connected has been explored in \cite{GrossReeder}. 
 We note that the case of $\SO_{2n + 1}$ has already been considered in \cite{Adrian}, and in particular \cite[Section 5]{Adrian} gives a special case of Proposition \ref{prop-decomp-pix}.

For the rest of the section, assume $\Gk$ is split, absolutely simple, and adjoint with root datum $(X = \Hom_k(S, \mathbb{G}_m), \Phi, \check X, \check \Phi)$. Choose a set of simple roots $\al_1, \dots, \al_\ell \in \Phi$, and let $\check\omega_1, \dots, \check\omega_\ell$ be the corresponding fundamental coweights. Since $\Gk$ is adjoint, the cocharacter group $\check X$ is generated by $\check\omega_1, \dots, \check\omega_\ell$. Let $\al_0 \in \Phi$ be the highest root with respect to our root basis, and define positive integers $c_1, \dots, c_\ell$ by $\al_0 = \sum_{i = 1}^\ell c_i\al_i$. We also set $c_0 = 1$. Equivalently, the $c_i$ are the labellings of the affine Dynkin diagram for $\Gk$.

\subsection{An invariant polynomial}

In this section we use Vinberg theory of graded Lie algebras to prove a lemma that will allow us to define an invariant polynomial for the action of $H_x$ on $\checkVx(\resk)$. For relevant background on Vinberg theory, see, e.g., \cite{Vinberg}, \cite{Panyushev}, \cite{RLYG}. For its relation to Moy--Prasad filtrations and its relevance in the current context, see, e.g., \cite[Section 4]{ReederYu}, \cite{FintzenMP}, \cite[Section 3.3]{RomanoThesis}.

In this section we let $L$ be any field, and let $\GZ$ be the split, absolutely simple adjoint group over $L$ whose root datum is the same as the root datum of $\Gk$. Let $\TZ$ be a maximal split torus in $\GZ$. We identify $X$ with the character group of $\TZ$, $\Phi$ with the roots of $\GZ$ with respect to $\TZ$, and make similar identifications for the other notation defined above. 
Given $\al = \sum_{i = 1}^\ell a_i\al_i \in \Phi$, let $\Ht(\al) = \sum_{i = 1}^\ell a_i$.
Choose a Chevalley basis for $\Lie(\GZ)$ with root vectors $\{E_\al \mid \al \in \Phi\}$, and 
let $\{u^L_\al: \mathbb{G}_a \to \GZ \mid \al \in \Phi\}$ be corresponding root-group homomorphisms.

Let $h = \sum_{i = 0}^\ell c_i$ be the Coxeter number of $\GZ$. Define subspaces of $\Lie(\GZ)$ by 
\begin{align*}
\Lie(\GZ)_0 &= \Lie(\TZ)\\
\Lie(\GZ)_j &= \spn_L\{E_\al \mid \Ht(\al) \equiv j \mod h\} \text{ if } j \in \{1, \dots, h -1\}.
\end{align*}
Then 
\begin{equation*}
\Lie(\GZ) = \underset{j \in \bZ/h\bZ}\bigoplus \Lie(\GZ)_j
\end{equation*}
is an $h$-graded Lie algebra. The Vinberg representation associated to this grading is the representation of $\TZ$ on $\VZ := \Lie(\GZ)_1$, where $\TZ$ acts by restriction of the adjoint action. 
Let $E_0 = E_{-\al_0}$, and for $i \in \{1, \dots, \ell\}$, let $E_i = E_{\al_i}$. Then $\{ E_0, E_1, \dots, E_\ell\}$ forms a basis for $\VZ$.
Let $\check\VZ$ be the dual representation, and let $\{F_0, \dots, F_\ell\}$ be the basis of $\check\VZ$ dual to $\{E_0, \dots, E_\ell\}$. We can explicitly write the action of $\TZ$ on this basis as follows. 

Every element of $\TZ(L)$ may be written as $\prod_{i = 1}^\ell \check\omega_i(t_i)$ for some $t_1, \dots, t_\ell \in L$. 
The action of such an element is given by
\begin{equation}\label{eqn-Gx-action}
\prod_{i = 1}^\ell \check\omega_i(t_i)\cdot \sum_{i = 0}^\ell a_iF_i = (\prod_{i = 1}^\ell t_i^{c_i})a_0F_0 + \sum_{i = 1}^\ell a_it_i^{-1}F_i.
\end{equation}

\begin{Lem}\label{lem-deltaZ}
Let $\Delta_L: \check\VZ \to L$ be defined by $\Delta_L(\sum_{i = 0}^\ell a_iF_i) = \prod_{i = 0}^\ell a_i^{c_i}$. Then 
\begin{enumerate}
\item[1.]
$\Delta_L$ is an invariant polynomial for the action of $\TZ$ on $\check\VZ$.
\item[2.] Suppose $\lam, \lam' \in \check\VZ$ with $\Delta_L(\lam) = \Delta_L(\lam') \neq 0$. Then $\lam$ and $\lam'$ are in the same $\TZ$-orbit.
\end{enumerate}
\end{Lem}

\textit{Proof}. Using (\ref{eqn-Gx-action}), it's not hard to see that $\Delta_L$ is invariant under the action of $\TZ$.
Let $\lam = \sum a_iF_i$ and $\lam' = \sum b_iF_i$, and suppose $\Delta_L(\lam) = \Delta_L(\lam') \neq 0$. Let $\mu = \prod \check{\omega}_i(a_ib_i^{-1})$. Then 
\begin{align*}
\mu\cdot a_iF_i &= b_iF_i,\\
 \mu \cdot a_iF_0 &= a_0\prod_{i = 1}^\ell a_i^{c_i}b_i^{-c_i}F_0\\
 &= b_0f_0
 \end{align*}
 since $a_0\prod_{i = 1}^\ell a_i^{c_i} = b_0\prod_{i = 1}^\ell b_i^{c_i}$. Thus $\mu\cdot \lam = \lam'$, and $\lam$ and $\lam'$ are in the same $\TZ$-orbit.
 \qed

Let $\mathscr{N}$ be the subgroup of $\GZ$ generated by $\TZ$ and $\{u_\al^L(1)u^L_{-\al}(-1)u^L_\al(1) \mid \al \in \Phi\}$. Let $\MZ$ be the subgroup of $\mathscr{N}$ stabilizing $\VZ$ under the adjoint action, i.e. 
$\MZ = \{n \in \mathscr{N}$ such that $n\cdot v \in \VZ$ for all $v \in \VZ\}$. If we let $W$ be the Weyl group of $\GZ$, then $\MZ/\TZ$ is isomorphic to the subgroup $W_0$ of $W$ that preserves the subset $\{-\al_0, \al_1, \dots, \al_\ell\}$ of $\Phi$ (cf. \cite[Proposition 6.4.2 and Theorem 7.2.2]{Carter}). 
Clearly $\MZ$ acts on $\VZ$ by restriction of the adjoint action, and so $\check\VZ$ also forms a representation of $\MZ$.

\begin{Prop}\label{prop-deltaM}
The polynomial $\Delta_L$ is invariant for the action of $\MZ$ on $\check\VZ$. 
\end{Prop}

\textit{Proof}. We may choose coset representatives $\{m_{w} \mid w \in W_0\}$ of $\MZ/\TZ$ such that for every $\al \in \Phi$, there exists an integer $\ep_{w, \al} \in \{-1, 1\}$ such that
\begin{equation*}
\Ad(m_w)(E_\al) = \epsilon_{w, \al}E_{w(\al)}.
\end{equation*} 
Further, the integers $\epsilon_{w, \al}$ can be defined over $\bZ$ in a way that is independent of the field of definition $L$ (see, e.g.,\cite[Proposition 6.4.2]{Carter}).
Since we already know that $\Delta_L$ is invariant under $\TZ$, it suffices to show that $\Delta_L(m_w\cdot f) = \Delta_L(f)$ for all $w \in W_0, f \in \check\VZ$. 

For ease of exposition, let $\be_0 = -\al_0$, and for $j \in \{1, \dots, \ell\}$ let $\be_j = \al_j$. Then each $w \in W_0$ permutes the set $\{\be_0, \dots, \be_\ell\}$, thus
 determining a permutation of $\{0, 1, \dots, \ell\}$ . By abuse of notation we denote this permutation by $w$. 

Given $w \in W_0$, we have $(m_w)^{-1} = m_{w^{-1}}t_w$ for some $t_w \in \TZ$, and 
\begin{align*}
m_w \cdot F_j &= \ep_{w^{-1}, w(\be_j)}\be_{w(j)}(t_w)F_{w(j)}\\ 
&= t_w^{-1} \cdot (\ep_{w^{-1}, w(\be_j)}F_{w(j)})
\end{align*}
for all $j \in \{0, \dots, \ell\}$.
Since $\Delta_L$ is invariant under $\TZ$, we have
\begin{align*}
\Delta_L(m_w \cdot \sum_{j = 0}^\ell a_jF_j) &= \Delta_L(t_w^{-1} \cdot \sum_{j = 0}^\ell \ep_{w^{-1}, w(\be_j)}a_jF_{w(j)})\\
&= \Delta_L(\sum_{j = 0}^\ell \ep_{w^{-1}, w(\be_j)}a_jF_{w(j)})\\
&= \prod_{j = 0}^\ell (\varep_{w^{-1}, w(\be_j)})^{c_j}\Delta_L(\sum a_jF_{w(j)}).
\end{align*}

Since $w \in W_0$, we have that $c_j = c_{w(j)}$ for all $j$ (to see this, note that $w$ must act as a symmetry of the affine Dynkin diagram of $\GZ$), so $\Delta_L(\sum a_jF_{w(j)}) = \Delta_L(\sum a_jF_j)$. 
Thus it suffices to show that
\begin{equation}\label{eqn-prod}
\prod_{j = 0}^\ell (\ep_{w, \be_j})^{c_j} = 1
\end{equation}
for all $w \in W_0$.
As mentioned above, the coefficients $\ep_{w, \al}$ are defined over $\bZ$ in a way that does not depend on $L$, and thus it suffices to prove (\ref{eqn-prod}) under the assumption that 
$L$ is algebraically closed and of characteristic zero. 

For the last step of the proof, we use a result of Panyushev about invariant polynomials on $\VZ$. So we now consider $\Delta_0: \VZ \to L$ defined by $\Delta_0(\sum_{j = 0}^\ell a_iE_i) = \prod_{j = 0}^\ell a_i^{c_i}$. Similar to the above,  it's easy to check that $\Delta_0$ is an invariant polynomial for the action of $\TZ$ on $\VZ$, and that 
\begin{equation*}
\Delta_0(m_w\cdot v) = \prod_{j = 0}^\ell (\ep_{w, \be_j})^{c_j} \Delta_0(v)
\end{equation*}
for all $w \in W_0, v \in \VZ$. 
Since $\Delta_0$ is invariant for the action of $\TZ$, \cite[Theorem 3.5(i)]{Panyushev} tells us that $\Delta_0$ is the restriction to $\VZ$ of a polynomial on $\Lie(\GZ)$ that is invariant under the adjoint action of $\GZ$. In particular, $\Delta_0$ is invariant under the action of $\MZ$, and thus $\prod_{j = 0}^\ell (\ep_{w, \be_j})^{c_j} = 1$ for all $w \in W_0$. But as shown above, this tells us that $\Delta_L$ is invariant under the action of $\MZ$. \qed

\subsection{The action of $H_x$}\label{section-Hx}

We return to the notation at the beginning of Section \ref{section-ssc}, so $\Gk$ is a split, adjoint, simple group over $k$ with root system $\Phi$ and a fixed choice of simple roots $\al_1, \dots, \al_\ell$. Given $\al \in \Phi$, let $U_\al \subset \Gk$ be the corresponding root group, and choose root-group homomorphisms $u_\al: \mathbb{G}_a \to U_\al$. 
These choices yield a hyperspecial point $x_0 \in \cA(S, k)$, which allows us to identify $\cA(S, k)$ with $\check X \otimes \bR$.
Let $\check\rho = \sum_{i = 1}^\ell \check\omega_i$, and let $h$ be the Coxeter number of $\Gk$. 
In this section, we fix 
\begin{equation*}
x = x_0 + \frac{1}{h}\check\rho,
\end{equation*}
 which is the barycenter of the fundamental alcove $\mathcal{C}$ corresponding to our choice of simple roots. For this choice of $x$, the subgroup $\Gk(k)_{x, 0}$ is an Iwahori subgroup, $H_x$ is the normalizer in $\Gk(k)$ of this Iwahori subgroup, and $\Gx$ is a rank-$\ell$ torus defined over $\resk$.

We can write $J_x$ and $J_x^+$ explicitly as follows. Let $\cO$ be the ring of integers of $k$, and let $\varpi \in \cO$ be a uniformizer. Let $\Phi^+$ be the set of positive roots with respect to $\al_1, \dots, \al_\ell$. Let
\begin{equation*}
S(k)_1 = \{s \in S(k) \mid \chi(s) \in 1 + \varpi\cO \text{ for all } \chi \in X\}.
\end{equation*}
Then 
\begin{align*}
J_x =& \langle S(k)_1, u_\al(y), u_{-\al}(\varpi y) \mid \al \in \Phi^+, y \in \cO\rangle,\\
J_x^+ =& \langle S(k)_1, u_{\al_i}(\varpi y), u_\al(y), u_{-\be}(\varpi y), u_{-\al_0}(\varpi^2 y) \mid 1 \leq i \leq \ell,\\ & \al \in \Phi^+\setminus \{\al_1, \dots, \al_\ell\}, \be \in \Phi^+\setminus \{\al_0\}, y \in \cO \rangle.
\end{align*}

To understand the simple supercuspidal representations of $G(k)$, we must understand the orbits for the action of $H_x$ on $\checkVx(\resk)$. 
Given $i \in \{1, \dots, \ell\}$, let 
\begin{equation*}
e_i = u_{\al_i}(1)J_x^+ \in J_x/J_x^+,
\end{equation*} 
and let $e_0 = u_{-\al_0}(\varpi)J_x^+$. Then $\{e_0, e_1, \dots, e_\ell \}$ forms
a basis for $\Vx(\resk) = J_x/J_x^+$ over $\resk$. 
Let $\{f_0, \dots, f_\ell\}$ be the dual basis of $\checkVx(\resk)$.

The quotient $H_x/J_x$ forms the $\resk$-points a finite reductive group with identity component $\Gx(\resk)$.
The component group of $H_x/J_x$, which we may identify with $H_x/\Gk_{x, 0}(k)$, is isomorphic to the stabilizer $\Omega$ of $\mathcal{C}$ in the extended affine Weyl group $\widetilde{W}$ of $\Gk$. More explicitly, we let $M$ be the subgroup of $\Gk(k)$ given by
\begin{equation*}
M = \la S(k), u_\al(1)u_{-\al}(-1)u_{\al}(1) \mid \al \in \Phi \ra.
\end{equation*}
 Then there is a natural surjective map $M \to \widetilde{W}$ (see \cite[Section 2.1]{IwahoriMatsumoto}). For each $\sigma \in \Omega$, we choose a preimage $n_\sigma$ of $\sigma$ in $M$ under this map. Then $\{n_\sigma \mid \sigma \in \Omega\}$ gives a complete set of coset representatives for $H_x/\Gk_{x, 0}(k)$ (\cite[Corollary 2.19]{IwahoriMatsumoto}).

We next use the previous section to describe an invariant polynomial for the action of $H_x$ on $\checkVx(\resk)$.

\begin{Lem}\label{lem-delta}
Let $\Delta: \checkVx(\resk) \to \resk$ be defined by $\Delta(\sum_{i = 0}^\ell a_if_i) = \prod_{i = 0}^\ell a_i^{c_i}$. Then 
\begin{itemize}
\item[1.] $\Delta$ is an invariant polynomial for the action of $H_x$ on $\checkVx(\frakf)$ (and thus also for the action of $\Gx(\frakf)$). 
\item[2.] A vector $\lam \in \checkVx(\resk)$ is stable if and only if $\Delta(\lam) \neq 0$. 
\item[3.] Suppose $\lam, \lam' \in \checkVx(\resk)$ are stable. Then $\lam$ and $\lam'$ are in the same $\Gx(\resk)$-orbit if and only if $\Delta(\lam) = \Delta(\lam')$.
\end{itemize}
\end{Lem}

\textit{Proof}.  
The second statement follows from the discussion in \cite[Section 2.6]{ReederYu}.
For the other statements, let $\GZ, \TZ, \VZ$, and $\check\VZ$ be as defined in the previous section, with $L = \resk$. Then there is a group isomorphism $\TZ \to \Gx$ and a linear isomorphism $\phi: \check\VZ \to \checkVx(\resk)$ given by $\phi(F_i) = f_i$ for $i \in \{0, \dots, \ell\}$ (\cite[Proposition 3.12]{RomanoThesis}). Under these isomorphisms, the action of $\TZ(\resk)$ on $\check\VZ$ corresponds to the action of $\Gx(\resk)$ on $\checkVx(\resk)$. Thus $\Delta$ is an invariant polynomial for the action of $\Gx(\resk)$ by the first part of Lemma \ref{lem-deltaZ}.
To finish the proof of the first statement, we must show that $\Delta(n_\sigma \cdot f) = \Delta(f)$ for all $\sigma \in \Omega$. So fix $\sigma \in \Omega$, and let $w \in W$ be the image of $\sigma$ under the natural projection $\widetilde{W} \to W$. By \cite[Section 2.4]{IwahoriMatsumoto}, we have $n_\sigma \cdot\phi(F) = \phi(m_w\cdot F)$ for all $F \in \check\VZ$. Thus the first statement follows from Proposition \ref{prop-deltaM}.
Using the isomorphism $\phi$, the third statement follows from the second statement and the second part of Lemma \ref{lem-deltaZ}. 
 \qed

\begin{Prop}\label{prop-orbit-reps}
\begin{enumerate}
\item[1.] Suppose $\lam, \lam' \in \checkVx(\resk)$ are stable. Then $\lam$ and $\lam'$ are in the same $H_x$-orbit if and only if $\Delta(\lam) = \Delta(\lam')$.
\item[2.] Let $q$ be the size of the residue field $\resk$. Then the set $I(x)$ contains exactly $q - 1$ equivalence classes under $\sim$. If we fix a nontrivial character $\chi_0$ of $\resk^+$ and let $\lam_c = cf_0 + \sum_{i = 1}^\ell f_i$, then $\{ (\lam_c, \chi_0) \mid c \in \resk^\times\}$ is a complete set of representatives for $I(x)/\sim$. 
\end{enumerate}
\end{Prop}

\textit{Proof}. This follows directly from Corollary \ref{cor-equiv-classes} and Lemma \ref{lem-delta}, where we are using the fact that if two vectors in $\checkVx(\resk)$ are in the same $\Gx(\resk)$-orbit, then they are in the same $H_x$-orbit.\qed

Using the previous results, we can explicitly describe the elements of each equivalence class in $I(x)$:

\begin{Prop}\label{prop-ssc-classification}
Let $(\lam, \chi), (\lam', \chi') \in I(x)$, and choose $c \in \resk^\times$ such that $\chi' = c\chi$. Then $(\lam, \chi) \sim (\lam', \chi')$ if and only if $\Delta(\lam)\Delta(\lam')^{-1} = c^h$. 
\end{Prop}

\textit{Proof}. Let $\lam = \sum a_if_i$ and $\lam' = \sum b_if_i$. Suppose $(\lam, \chi) \sim (\lam', \chi')$. Then by Theorem \ref{thm-input} there exists $g \in H_x$ such that $g\lam = c\lam'$. Since $\Delta$ is homogeneous of degree $h$, we have $\Delta(c\lam') = c^h\Delta(\lam')$. Since $\Delta$ is $H_x$-invariant, we have 
\begin{equation*}
\Delta(\lam) = \Delta(g\lam) = c^h\Delta(\lam'),
\end{equation*}
 and $\Delta(\lam)\Delta(\lam')^{-1} = c^h$ as claimed.
Now suppose $\Delta(\lam)\Delta(\lam')^{-1} = c^h$. Then $\Delta(\lam) = \Delta(c\lam')$, so by Proposition \ref{prop-orbit-reps}, $\lam$ and $c\lam'$ are in the same $H_x$-orbit, and so by Theorem \ref{thm-input} we have $(\lam, \chi) \sim (\lam', \chi')$.\qed

\subsection{Decomposing $\pi_x(\lam, \chi)$}\label{section-decomp-ssc}

We continue with the notation of the previous section, so $x = x_0 + \frac{1}{h}\check\rho$ is the barycenter of the fundamental alcove corresponding to our choice of simple roots. In addition, for the rest of the section, we fix a stable vector $\lam \in \checkVx(\resk)$ and a nontrivial character $\chi$ of $\resk^+$.
Our next goal is to describe the output of the construction of \cite{ReederYu}, i.e. the irreducible subrepresentations of $\pi_x(\lam, \chi)$. To do so, we must first understand the subgroup $H_\lam = \Stab_{H_x}(\lam)$. 

Note that the following lemma can be deduced from \cite[Proposition 5.5.1]{Rostami}, but here we give a comparatively simple proof using the invariant polynomial $\Delta$. As in Section \ref{section-Hx}, we denote by $\Omega$ the stabilizer of the fundamental alcove $\mathcal{C}$ in the extended affine Weyl group $\widetilde{W}$.

\begin{Lem}\label{lem-iso-omega}
\begin{enumerate}
\item[1.] The natural inclusion $H_\lam \hookrightarrow H_x$ induces an isomorphism $H_\lam/J_x \simeq H_x/\Gk(k)_{x, 0}$. In particular, $H_\lam/J_x \simeq \Omega$. 
\item[2.] We can choose coset representatives for $H_\lam/J_x$ in $M$.
\end{enumerate}
\end{Lem}

\textit{Proof}. 
Let $\lam = \sum_{i = 0}^\ell a_if_i$, and let $\iota: H_\lam/J_x \to H_x/\Gk(k)_{x, 0}$ be the homomorphism induced by the inclusion map. 
Note that $\Stab_{\Gx(\resk)}(\lam)$ is trivial: this follows from the fact that $\Gx$ acts as in (\ref{eqn-Gx-action}) and the fact that the coefficients $a_i$ are all nonzero. Using this,
we have  that
$H_\lam \cap \Gk(k)_{x, 0} \subset J_x$, which implies that $\iota$ is injective. 
To show $\iota$ is surjective, let $h \in H_x$. 
Then by Lemma \ref{lem-delta} part 1 we have $\Delta(h^{-1}\lam) = \Delta(\lam)$, 
so by Lemma \ref{lem-delta} part 3, there exists $g \in \Gk(k)_{x, 0}$ such that $h^{-1}\cdot \lam = g\lam$. Thus $h g \in H_\lam$, and the map $H_\lam/J_x \to H_x/\Gk(k)_{x, 0}$ is an isomorphism. 

As in Section \ref{section-Hx}, we can choose coset representatives $\{n_\sigma \mid \sigma \in \Omega\}$ of $H_x/\Gk(k)_{x, 0}$ in $M$. By the previous paragraph, a set of coset representatives for $H_\lam/J_x$ is then given by $\{n_\sigma g_\sigma\}$ where $g_\sigma$ is any element of $\Gk(k)_{x, 0}$ such that $n_\sigma^{-1}\lam = g_\sigma\lam$. Since the action of $\Gk(k)_{x, 0}$ factors through $\Gk(k)_{x, 0}/J_x = \Gx(\resk)$ and $M \cap \Gk(k)_{x, 0}$ maps surjectively onto $\Gx(\resk)$, we can choose each $g_\sigma$ to be in $M$.
\qed

\begin{Lem}\label{lem-semidirect}
$H_\lam/J_x^+$ is isomorphic to a semidirect product $H_\lam/J_x \ltimes J_x/J_x^+$.
\end{Lem}

\textit{Proof}. 
Let $\{h_\sigma \mid \sigma \in \Omega\}$ be a set of coset representatives for $H_\lam/J_x$ with $h_\sigma \in M$ for all $\sigma$. 
We claim the map 
\begin{equation*}
H_\lam/J_x \to H_\lam/J_x^+
\end{equation*}
 given by $h_\sigma J_x \mapsto h_\sigma J_x^+$
  is a homomorphism. Indeed, if $\sigma, \tau \in \Omega$, then $h_\sigma h_\tau (h_{\sigma\tau})^{-1} \in M \cap J_x$ since $H_\lam/J_x \simeq \Omega$ is abelian.
But $M \cap J_x \subset J_x^+$, so $(h_\sigma h_\tau) J_x^+ = h_{\sigma\tau}J_x^+$. Thus the short exact sequence
\begin{equation*}
1 \to J_x/J_x^+ \to H_\lam/J_x^+ \to H_\lam/J_x \to 1
\end{equation*}
splits, which proves the lemma.
\qed

\begin{Prop}\label{prop-decomp-pix}
The representation $\ind_{J_x}^{H_\lam}\chi_\lam$ decomposes as a direct sum of $\lvert \Omega \rvert$ characters. 
In particular, $\Irr(\lam, \chi)$ contains $\lvert \Omega \rvert$ irreducible supercuspidal representations: the representations of the form $\ind_{H_\lam}^{G(k)} \tilde\chi$ for the characters $\tilde\chi$ of $H_\lam$ such that $\tilde\chi\mid_{J_x} = \chi_\lam$.
\end{Prop}

\textit{Proof}. 
First note that there is a natural bijection between characters of $H_\lam/J_x \simeq \Omega$ and the characters $\tilde\chi$ of $H_\lam$ such that $\tilde\chi\mid_{J_x} = \chi_\lam$. Indeed, every such character $\tilde\chi$ of $H_\lam$ factors through $H_\lam/J_x^+$, so using the decomposition of Lemma \ref{lem-semidirect}, we obtain a character of $H_\lam/J_x$ by restriction. And if $\rho$ is a character of $H_\lam/J_x$, then it is an easy exercise to show that $\rho$ extends uniquely to a character of $H_\lam/J_x^+ \simeq H_\lam/J_x \ltimes J_x/J_x^+$ whose restriction to $J_x/J_x^+$ is $\chi_\lam$ (here we use the fact that $\chi_\lam^h = \chi_\lam$ for all $h \in H_\lam$, similar to \cite[Section 8.2]{Serre}). Since $\Omega$ is abelian, this shows that there are exactly $\lvert \Omega \rvert$ characters $\tilde\chi$ of $H_\lam$ such that $\tilde\chi\mid_{J_x} = \chi_\lam$.

Now fix a character $\tilde\chi$ of $H_\lam$ as above. Since $H_\lam/J_x$ is finite, 
we can use Frobenius reciprocity to conclude that
\begin{equation*}
\Hom_{H_\lam}(\tilde\chi, \ind_{J_x}^{H_\lam} \chi_\lam) \simeq \Hom_{J_x}(\chi_\lam, \chi_\lam)
\end{equation*}
(see, e.g., \cite[Section 2.4]{BushnellHenniartGL2}). 
Thus $\tilde\chi$ appears as a subrepresentation of $\ind_{J_x}^{H_\lam}\chi_\lam$. Since $\ind_{J_x}^{H_\lam}\chi_\lam$ has dimension $[H_\lam:J_x]$, every irreducible subrepresentation is of this form. The result then follows from \cite[Proposition 2.4]{ReederYu}.
 \qed

\begin{remark}\label{rem-psi-notation}
Note that given $(\lam, \chi) \in I(x)$, the proof of Proposition \ref{prop-decomp-pix} associates to each character $\psi$ of $\Omega$ a subrepresentation of $\pi_x(\lam, \chi)$. We denote this subrepresentation by $\pi_x(\lam, \chi, \psi)$.
\end{remark}

\begin{Cor}\label{cor-formal-deg}
For $x = x_0 + \frac{1}{h}\check\rho$, the construction of \cite{ReederYu} gives a total of $(q - 1)\lvert \Omega \rvert$ supercuspidal representations. With respect to a Haar measure $\mu$ on $\Gk(k)$, each of these representations has 
 formal degree $\frac{1}{\lvert \Omega \rvert \mu(J_x)}$. 
 \end{Cor}
 
 \textit{Proof}. The first statement follows from combining Propositions \ref{prop-orbit-reps} and \ref{prop-decomp-pix}. The second statement follows from \cite[Proposition 2.4]{ReederYu}, where we are using the fact that for all $(\lam, \chi) \in I(x)$, each irreducible subrepresentation of $\ind_{J_x}^{H_\lam}\chi_\lam$ has dimension one, and $[H_\lam: J_x] = \lvert \Omega \rvert$. \qed

\subsection{Unramified extensions}\label{section-unram}

In the last section of the paper, we will consider Langlands parameters associated to the representations of Proposition \ref{prop-decomp-pix} following \cite{GrossReeder}. In order to do this, we now relate the simple supercuspidals of Proposition \ref{prop-decomp-pix} to certain simple supercuspidals for $\Gk(k')$, where $k'$ is a finite unramified extension of $k$. 

For every positive integer $m$, let $k_m$ be the degree-$m$ unramified extension of $k$, and let $\resk_m$ be the degree-$m$ extension of $\resk$. 
The base change $\Gk_{k_m}$ is a split reductive group over $k_m$ with maximal torus $S_{k_m}$. 
Because $\Gk$ is split, we may identify the apartment $\cA(S_{k_m}, k_m)$ with the apartment $\cA(S, k)$, and thus we may think of $x = x_0 + \frac{1}{h}\check\rho$ as lying in $\cA(S_{k_m}, k_m)$.  
For every positive integer $m$, we let $J_m = \Gk(k_m)_{x, r(x)}$.
If we consider Moy--Prasad subgroups in $\Gk(k_m)$, we have 
\begin{equation*}
\Gk(k_m)_{x, 0}/\Gk(k_m)_{x, r(x)} \simeq \Gx(\resk_m) \text{ and } J_{m}/J_m^+ \simeq \Vx(\resk_m) \simeq \Vx(\resk) \otimes \resk_m.
\end{equation*}

Thus for each $m$, we can use the construction of \cite{ReederYu} with respect to our choice of $x$ to build representations of $\Gk(k_m)$.
We denote the corresponding set of the input as $I_m(x)$, so $I_m(x)$ is the set of pairs $(\lam, \chi)$ such that $\lam$ is a stable vector in $\checkVx(\resk_m)$ and $\chi$ is a nontrivial character of $\resk_m^+$. 

Fix $(\lam, \chi) \in I(x) = I_1(x)$. 
Note that we may consider $\lam$ as an element of $\checkVx(\resk_m) = \checkVx(\resk) \otimes \resk_m$, and by the definition of stable, $\lam$ is stable as an element of this space. Let $\text{Tr}_m: \resk_m \to \resk$ denote the trace map, which is a homomorphism of additive groups. Since $\text{Tr}_m$ is surjective, we see that
\begin{equation*}
\chi_m := \chi \circ \text{Tr}_m: \resk_m \to \bC^\times
\end{equation*}
 is a nontrivial character of $\resk_m$, so $(\lam, \chi_m) \in I_m(x)$.

Further, we can define a canonical bijection between the set of representations $\Irr(\lam, \chi)$ of $\Gk(k)$ and the set of representations $\Irr(\lam, \chi_m)$ of $\Gk(k_m)$ as follows: let $H_{m} = \Stab_{\Gk(k_m)}(x)$, and let $H_{\lam, m} = \Stab_{H_m}(\lam)$. Note that $\Omega$ is determined by the root datum of $\Gk$ and does not depend on $m$. By Lemma \ref{lem-iso-omega}, the quotient $H_{\lam, m}/J_m$ is canonically isomorphic to $\Omega$ and so does not depend on $m$. Using Proposition \ref{prop-decomp-pix} we have a bijection between $\Irr(\lam, \chi)$ and $\Irr(\lam, \chi_m)$ that, in the notation of Remark \ref{rem-psi-notation}, is given by
\begin{equation*}
\pi_x(\lam, \chi, \psi) \mapsto \pi_x(\lam, \chi_m, \psi).
\end{equation*}
In Section \ref{section-LLC} we will make an assumption about the Langlands parameters of these representations analogous to an assumption made in \cite{GrossReeder}.

\begin{remark}
The use of the trace map in this section is analogous to the use of $\tau_m$ in \cite[Section 9.5]{GrossReeder}.
\end{remark}

\section{Langlands parameters for simple supercuspidal representations of adjoint groups}\label{section-LPs}

We finish the paper by describing properties of the Langlands parameters for the simple supercuspidal representations classified in Section \ref{section-ssc}.

\subsection{Langlands parameters and L-functions}\label{section-LP-defs}

We now review some background and set up notation following \cite{GrossReeder}, which the reader can refer to for more details. 
In this section, we assume that $\Gk$ is absolutely simple and split over $k$. 
Let $G^\vee$ be the ($\bC$-points of) the dual group of $\Gk$, and let $Z^\vee$ be the center of $G^\vee$. Let $\Weil$ be the Weil group of $k$, and let $\Inert \subset \Weil$ be the inertia subgroup of $\Weil$. Fix a geometric Frobenius $\Frob \in \Weil$. Since $\Gk$ is split, we may define a Langlands parameter for $\Gk$ to be a homomorphism
\begin{equation*}
\varphi: \Weil \times \SL_2(\bC) \to G^\vee
\end{equation*}
such that 
$\varphi$ is continuous on $\Inert$ and algebraic on $\SL_2(\bC)$, and such that $\varphi(\Frob)$ is semisimple. 
Two Langlands parameters are said to be equivalent if they are conjugate by an element of $G^\vee$. 
We write $\princ$ for the principal parameter, which has the following properties: $\princ\vert_{\Weil}$ is trivial, and $\princ(\begin{psmallmatrix} 1 & 1\\0 & 1\end{psmallmatrix})$ is regular unipotent in $G^\vee$.

A Langlands parameter $\varphi$ is called \textit{discrete} if the centralizer $A_\varphi := Z_{G^\vee}(\im \varphi)$ is finite. Given a discrete Langlands parameter $\varphi$, we think of an \textit{enhancement} of $\varphi$ as an irreducible representation of $A_\varphi$. Further, since $\Gk$ is split, an enhancement $\rho$ is said to be \textit{relevant} for $\Gk$ if $\rho\vert_{Z^\vee}$ is trivial. We write $\Phi(G)$ for the equivalence classes of enhanced discrete Langlands parameters $(\varphi, \rho)$ that are relevant for $\Gk$. 

\begin{Lem}\label{lem-triv-rho}
Suppose $(\varphi, \rho) \in \Phi(\Gk)$ satisfies
$\dim \rho = \lvert A_\varphi/Z^\vee \rvert$. Then $A_\varphi = Z^\vee$, and $\rho$ is the trivial representation.
\end{Lem}

\textit{Proof}. By definition, $\rho$ is an irreducible representation of $A_\varphi$ trivial on $Z^\vee$, so we may think of it as an irreducible representation of $A_\varphi/Z^\vee$. Since every irreducible representation of $A_\varphi/Z^\vee$ is a subrepresentation of the regular representation, $\rho$ has dimension less than or equal to $\lvert A_\varphi/Z^\vee \rvert$, with equality if and only if $A_\varphi/Z^\vee$ is the trivial group. The result follows.
\qed

Fix a discrete Langlands parameter $\varphi$. Let $\frakg^\vee = \Lie(G^\vee)$, and write $\Ad$ for the adjoint representation of $G^\vee$. The parameter $\varphi$ determines a representation of $\Weil$ on $\frakg^\vee$ via $\Ad \circ \varphi$. Let $(\frakg^\vee)^{\Inert}$ be the fixed points for the action of $\Inert$ on $\frakg^\vee$.
Similarly $\varphi$ determines a representation of $\SL_2(\bC)$ on $\frakg^\vee$. 
Let $U$ be the subspace of $(\frakg^\vee)^{\Inert}$ fixed by the action of $\begin{psmallmatrix}1 & 1\\0 & 1\end{psmallmatrix} \in \SL_2(\bC)$. Let 
\begin{equation*}
F_q = \Ad \circ \varphi((\Frob, \begin{psmallmatrix} q^{-1/2} & 0\\0 & q^{1/2} \end{psmallmatrix})).
\end{equation*}
We define the adjoint L-function $L(\varphi, s)$ be
\begin{equation*}
L(\varphi, s) = \det(I - q^{-s}F_q\vert_U)^{-1},
\end{equation*}
where $I$ is the identity automorphism of $U$. 
It is helpful to decompose $L(\varphi, s)$ as a product in the following way. We may decompose $\frakg^\vee$ as
\begin{equation*}
\frakg^\vee = \underset{n \geq 0}\bigoplus~ \frakg^\vee_n \otimes \Sym^n,
\end{equation*}
where $\Sym^n = \Sym^n(\bC^2)$ is the irreducible $(n + 1)$-dimensional representation of $\SL_2(\bC)$ and $\frakg^\vee_n$ is a semisimple complex representation of $\Weil$. 
Let 
\begin{equation*}
U_n = (\frakg^\vee_n)^{\Inert}.
\end{equation*}
Then we may rewrite $L(\varphi, s)$ as
\begin{equation*}
L(\varphi, s) = \prod_{n \geq 0} \det(I_n - q^{-\frac{n}{2} - s}\varphi(\Frob)\vert_{U_n})^{-1},
\end{equation*}
where $I_n$ is the identity automorphism of $U_n$.

As in \cite[Proposition 4.3]{GrossReeder}, for each $n \geq 0$ define $P_n(t)$ by
\begin{equation*}
P_n(t) = \det(I_n - t\varphi(\Frob)\vert_{U_n}),
\end{equation*} 
so that 
\begin{equation}\label{eqn-LPn}
L(\varphi, s) = \prod_{n \geq 0} P_n(q^{-s - n/2})^{-1}.
\end{equation}
Then $P_n(t) \in \bZ[t]$ for all $n$. Further, for $n \geq 1$, $P_n(t)$ may be factored as $\prod R_i(t^{e_i})$, where each $R_i \in \bZ[t]$, and $ne_i$ is even for all $i $ by \cite[Proposition 4.3]{GrossReeder}. (Note that \cite{GrossReeder} assumes that $k$ has characteristic zero, but the proof of Proposition 4.3 reduces to analysis of torsion automorphisms on complex semisimple Lie algebras, and in particular does not use the characteristic-zero assumption.)

\begin{Lem}\label{lem-P2}
If $P_2(t)$ is constant, then $\varphi(\SL_2(\bC)) = 1$.
\end{Lem}

\textit{Proof}. We follow the same logic as the proof of  \cite[Lemma 5.3]{GrossReeder}. Suppose $\varphi(\SL_2(\bC)) \neq 1$. Then $\varphi(\SL_2(\bC))$ is a subgroup of $G^\vee$ of type $A_1$, and $\Ad\vert_{\varphi(\SL_2(\bC))}$ is its adjoint representation. Since $\varphi(\Inert)$ centralizes $\varphi(\SL_2(\bC))$, this tells us that $(\frakg_2^\vee)^{\Inert} \neq 0$, and $P_2(t)$ is nonconstant. The result follows. \qed

As \cite[Section 2.2]{GrossReeder}, let
\begin{align*}
\al(\varphi) &= a(\varphi) + \underset{n \geq 1}\sum n\dim U_n,\\
\omega(\varphi) &= \underset{n \geq 0}\prod w(\frakg^\vee_n)^{n + 1}\underset{n \geq 1}\prod \det(-\varphi(\Frob)\vert_{U_n})^n
\end{align*}
where $a(\varphi)$ is the Artin conductor for the representation of $\Weil$ on $\frakg^\vee$, and $w(\frakg^\vee_n)$ is the local constant for the action of $\Weil$ on $\frakg^\vee_n$. 
We also write $b(\varphi)$ for the Swan conductor for the action of $\Weil$ on $\frakg^\vee$ and note that 
\begin{equation*}
a(\varphi) = \dim (\frakg^\vee/(\frakg^\vee)^\Inert) + b(\varphi).
\end{equation*}

\begin{Lem}\label{lem-Un}
Suppose $\varphi(\Weil) = 1$. Then
\begin{equation*}
\sum_{n \geq 1} n\dim U_n = \dim \frakg^\vee - \dim U.
\end{equation*}
\end{Lem}

\textit{Proof}. This is easily verified using the fact that $\dim \Sym^n = n + 1$, and the fact that $\dim U_n$ is the multiplicity of $\Sym^n$ as a subrepresentation of $\frakg^\vee$ under the action of $\varphi(\SL_2(\bC))$. So if $V_1, \dots, V_m$ are irreducible representations of $\SL_2(\bC)$ with $\frakg^\vee \simeq \oplus_{i = 1}^m V_i$, then 
\begin{equation*}
\sum_{n \geq 1} n\dim U_n = \sum_{i = 1}^m (\dim V_i -1).
\end{equation*} \qed

Let
\begin{equation*}
\varep(\varphi, s) = \omega(\varphi)q^{\al(\varphi)(\frac{1}{2} - s)},
\end{equation*}
and let
\begin{equation*}
\gamma(\varphi) = \frac{L(\varphi, 1)\varep(\varphi, 0)}{L(\varphi, 0)}
\end{equation*}
be the adjoint gamma value associated to $\varphi$ (cf. \cite[Section 4]{GrossReeder}).

For future use we consider the case of the principal parameter. Suppose $\Gk$ is adjoint of rank $\ell$, as in Section \ref{section-ssc}. Let $\princ$ be the principal parameter, and let $\mu_G$ be the absolute value of the Euler--Poincar\'e measure on $\Gk(k)$ (cf. \cite{Serre71}, \cite[Section 7.1]{GrossReeder}). 
By Lemma \ref{lem-Un}, we have
\begin{equation*}
\al(\princ) = \dim \frakg^\vee - \ell.
\end{equation*}
By \cite[Corollary 8.7]{Kostant}, we have
\begin{equation}
\frac{L(\princ, 1)}{L(\princ, 0)} = q^\ell \prod_{i = 1}^\ell \frac{q^{m_i} - 1}{q^{m_i + 1} - 1},
\end{equation}
where $m_1, \dots, m_\ell$ are the exponents of $G^\vee$.
Using \cite[Th\'eor\`eme 7]{Serre71}, this gives us
\begin{equation*}
\frac{L(\princ, 1)}{L(\princ, 0)} = \lvert \Omega\rvert q^\ell\mu_G(J_x),
\end{equation*}
where $J_x$ is as defined in Section \ref{section-Hx}.
Putting this together with the definition of $\gamma$, we have
\begin{equation}\label{eqn-gamma-princ}
\gamma(\princ) = \omega(\princ)\lvert \Omega \rvert q^{\frac{1}{2}(\dim \frakg^\vee + \ell)}\mu_G(J_x).
\end{equation}

\subsection{The local Langlands correspondence}\label{section-LLC}

For the rest of the paper we continue with the notation of Section \ref{section-ssc}, so $\Gk$ is split, absolutely simple, and adjoint. We now describe properties of the Langlands parameters associated to the representations $\pi_x(\lam, \chi, \psi)$ defined in Section \ref{section-decomp-ssc}. We note that the results of this section generalize results of \cite{Adrian} for $\SO_{2n + 1}$ in the case when the residue characteristic is sufficiently large.

As in Section \ref{section-unram}, given a positive integer $m$, we let $k_m$ be the degree-$m$ unramified extension of $k$. Let $\Weil_m$ be the Weil group of $k_m$. Given a discrete Langlands parameter $\varphi$, let $\varphi^m = \varphi\mid_{\Weil_m}$. Then $\varphi^m$ is a Langlands parameter for the base change $\Gk_{k_m}$ such that the actions of $\SL_2(\bC)$ and of $\Inert \subset \Weil_m$ are the same as the corresponding actions for $\varphi$. Thus we have
\begin{equation}\label{eqn-Lphi_m}
L(\varphi^m, s) = \prod_{n \geq 0} \det(I_n - (q^m)^{-\frac{n}{2} - s}\varphi(\Frob)^m\vert_{U_n})^{-1}.
\end{equation}
Since $\varphi$ is discrete, by \cite[Lemma 3.1]{GrossReeder}, we know that $\varphi(\Frob)$ has finite order, say $d$. Note that if $m \equiv 1 \mod d$, then 
$\im \varphi = \im \varphi^m$, so $A_\varphi = A_{\varphi^m}$ and an enhancement of $\varphi$ is the same as an enhancement of $\varphi^m$.

For each positive integer $m$, let $\Rep(\Gk({k_m}))$ be the set of (isomorphism classes of) irreducible discrete-series representations of $\Gk(k_m)$, and let $\Phi(G_{k_m})$ be the enhanced discrete Langlands parameters for $\Gk_{k_m}$ as defined in the previous section. 
For each such $m$, we assume
\begin{equation*}
\LL: \Rep(\Gk({k_m})) \to \Phi(G_{k_m})
\end{equation*} 
is a map satisfying the expected properties of the local Langlands correspondence for $\Gk_{k_m}$. We write $\LLG$ for $\LLk$. 
More specifically, to prove the results of this section, it suffices to make the following assumptions:
\begin{enumerate}
\item[A1.] Given a discrete-series representation $\pi \in \Rep(\Gk({k_m}))$, let $\Deg(\pi)$ be the formal degree of $\pi$ with respect to $\mu_{G_{k_m}}$, and let $\princ$ be the principal parameter for $\Gk_{k_m}$. If $\LL(\pi) = (\varphi, \rho)$, then 
\begin{equation}
\Deg(\pi) = \frac{\omega(\varphi)}{\omega(\princ)}\frac{\dim \rho}{\lvert A_\varphi/Z^\vee \rvert} \frac{\gamma(\varphi)}{\gamma(\princ)}.
\end{equation}
\item[A2.] With notation as in Section \ref{section-unram}, if $\LLG(\pi_x(\lam, \chi, \psi)) = (\varphi, \rho)$, then 
\begin{equation*}
\LL(\pi_x(\lam, \chi_m, \psi)) = (\varphi^m, \rho)
\end{equation*} 
for all $m \equiv 1 \mod 4\lvert \varphi(\Weil)\rvert$. 
\end{enumerate}
To justify these assumptions, note that 
A1 is the formal degree conjecture of \cite{HiragaIchinoIkeda} as reformulated in \cite{GrossReeder}, and A2 is an analogue of the assumptions (75) and (76) in \cite{GrossReeder}.

\begin{Thm}\label{thm-LP}
Let $(\varphi, \rho) = \LLG(\pi_x(\lam, \chi, \psi))$. Then 
\begin{enumerate}
\item[1.]
$\al(\varphi) = \dim \frakg^\vee + \ell$. 
\item[2.]  $\varphi(\SL_2(\bC)) = 1$.
\item[3.] The adjoint L-function $L(\varphi, s)$ is trivial. 
\item[4.] $A_\varphi = Z^\vee$, and $\rho$ is the trivial representation.
\end{enumerate}
\end{Thm}

\textit{Proof}. 
Let $d = \lvert \varphi(\Weil) \rvert$. 
Putting together assumptions A1 and A2, we have
\begin{equation}\label{eqn-FDC}
 \Deg(\pi_x(\lam, \chi_m, \psi)) =  \frac{\omega(\varphi^m)}{\omega(\princ^m)}\frac{\dim \rho}{\lvert A_{\varphi}/Z^\vee\rvert} \frac{\gamma(\varphi^m)}{\gamma(\princ^m)}
\end{equation}
for all $m \equiv 1 \mod 4d$. 

For all $m \geq 1$, let $\mu_m$ be the absolute value of the Euler--Poincar\'e measure on $\Gk(k_m)$. By Corollary \ref{cor-formal-deg}, the lefthand side of (\ref{eqn-FDC}) is $\frac{1}{\lvert \Omega \rvert \mu_m(J_m)}$. By (\ref{eqn-gamma-princ}), the righthand side is
\begin{equation*}
\left(\frac{\omega(\varphi^m)}{\omega(\princ^m)}\right)^2 \frac{\dim \rho}{\lvert A_\varphi/Z^\vee\rvert}\frac{1}{\lvert \Omega \rvert \mu_m(J_m)}q^{\frac{m}{2}(\al(\varphi^m) - \dim \frakg^\vee - \ell)}\frac{L(\varphi^m, 1)}{L(\varphi^m, 0)}.
\end{equation*}
Note that $\al(\varphi^m) = \al(\varphi)$ for all $m$, since $\al(\varphi)$ only depends on the action of $\Inert$. Also note that $\frac{\omega(\varphi^m)}{\omega(\princ^m)} \in \{1, -1\}$. Thus (\ref{eqn-gamma-princ}) implies that
\begin{equation}\label{eqn-Lratio}
\frac{L(\varphi^m, 1)}{L(\varphi^m, 0)} = \frac{\lvert A_\varphi/Z^\vee\rvert}{\dim \rho} q^{\frac{m}{2}(\dim \frakg^\vee + \ell - \al(\varphi))}
\end{equation}
for all $m \equiv 1 \mod 4d$. 

With notation as in Section \ref{section-LP-defs}, let 
\begin{align*}
f(t) &= \prod_{n \geq 0} P_n(t^{n/2}) \in \bZ[t] \text{ and}\\ 
g(t) &= \prod_{n \geq 0} P_n(t^{1 + n/2}) \in \bZ[t].
\end{align*} 
Note that these products are finite: $\frakg^\vee$ is finite dimensional, so the subspaces $\frakg_n^\vee$ are trivial for all but finitely many $n$. The fact that $f(t), g(t) \in \bZ[t]$ follows from \cite[Proposition 4.3]{GrossReeder}, as mentioned in the previous section.

Using (\ref{eqn-Lphi_m}), if $m \equiv 1 \mod d$, then 
\begin{align*}
L(\varphi_m, s) &= \prod_{n \geq 0} \det(I_n - (q^m)^{-\frac{n}{2} - s}\varphi(\Frob)\vert_{U_n})^{-1}\\
&= \prod_{n \geq 0} P_n((q^m)^{-\frac{n}{2} - s})^{-1},
\end{align*}
so $L(\varphi^m, 0) = f(q^{-m})^{-1}$ and $L(\varphi^m, 1) = g(q^{-m})^{-1}$. 
By (\ref{eqn-Lratio})
\begin{equation*}
f(q^{-m}) = \frac{\lvert A_\varphi/Z^\vee\rvert}{\dim \rho} q^{\frac{m}{2}(\dim \frakg^\vee + \ell - \al(\varphi))}g(q^{-m})
\end{equation*}
for infinitely many values of $m$, so 
\begin{equation*}
f(q^{-m})^2 = \left(\frac{\lvert A_\varphi/Z^\vee\rvert}{\dim \rho}\right)^2 (q^{-m})^{(\al(\varphi) - \dim \frakg^\vee - \ell)}g(q^{-m})^2
\end{equation*}
for infinitely many values of $m$. We may conclude that
\begin{equation}\label{eqn-equalpoly}
f(t)^2 = \left(\frac{\lvert A_\varphi/Z^\vee\rvert}{\dim \rho}\right)^2 t^{(\al(\varphi) - \dim \frakg^\vee - \ell)}g(t)^2
\end{equation}
as rational functions. This implies $\dim \frakg^\vee + \ell - \al(\varphi) = 0$, since $f(0) = 1$. The first part of the theorem follows.

Now note that $\deg g \geq \deg f$, with equality if and only if $P_n$ is constant for all $n \geq 0$. 
Thus (\ref{eqn-equalpoly}) implies that $P_n$ is constant for all $n \geq 0$. Since $P_n(0) = 1$, we have $P_n(t) = 1$ for all $n \geq 0$. By Lemma \ref{lem-P2}, this implies the second part of the theorem, and (\ref{eqn-LPn}) implies the third part.
Going back to (\ref{eqn-Lratio}), we see that $\dim \rho = \lvert A_\varphi/Z^\vee\rvert$, which implies the fourth part of the theorem by Lemma \ref{lem-triv-rho}.
\qed

\begin{remark}\label{rem-packet}
One of the conditions of the local Langlands correspondence says that the L-packet in $\Rep(G(k))$ corrresponding to a discrete Langlands parameter $\varphi$ should be in bijection with the irreducible representations of $A_\varphi/Z^\vee$. Thus a corollary of the fourth statement in Theorem \ref{thm-LP} is that each representation $\pi_x(\lam, \chi, \psi)$ lies in a singleton L-packet. This is in contrast with simple supercuspidals for simply connected $\Gk$, where the size of the L-packet depends on the residue cardinality $q$ (for an example see Theorem 1.1 and Remark 1.2 of \cite{HenniartOi}).
\end{remark}

\begin{Cor}\label{cor-swan}
Let $(\varphi, \rho) = \LLG(\pi_x(\lam, \chi, \psi))$. Then 
$(\frakg^\vee)^{\Inert} = 0$ and $b(\varphi) = \ell$. In particular $\varphi$ has minimal Swan conductor, in the language of \cite{GrossReeder}.
\end{Cor}

\textit{Proof}. The fact that $(\frakg^\vee)^\Inert = 0$ follows as in the proof of \cite[Lemma 5.3]{GrossReeder}. Since $U = 0$, we have $\al(\varphi) = a(\varphi)$, and since $(\frakg^\vee)^\Inert = 0$, we have 
\begin{equation*}
a(\varphi) = \dim \frakg^\vee + b(\varphi).
\end{equation*}
The fact that $b(\varphi) = 0$ then follows from the first part of Theorem \ref{thm-LP}. The fact that $\varphi$ has minimal Swan conductor then follows directly from the definition in \cite[Section 5]{GrossReeder}.
\qed

\begin{remark}
It is natural to wonder whether the methods of Theorem \ref{thm-LP} (in particular, the use of the polynomials $P_n(t)$ defined in \cite{GrossReeder}) can be used to prove similar results about Reeder--Yu representations corresponding to other choices of $x$. 
To do this, one would have to formulate an analogue of property A2 above. In particular, in the current setting, we used the fact that the irreducible subrepresentations of $\pi_x(\lam, \chi_m)$ are in canonical bijection with the irreducible representations of $\Omega$, no matter what $m$ is (recall that $\Omega$ depends only on the root datum of $\Gk$). For general $x$, the stabilizer of $\lam$ can change depending on $m$, so the situation is more complicated.
\end{remark}

\bibliographystyle{plain}
\bibliography{on-reeder-yu-const}

\end{document}